\documentclass[12pt]{amsart}

\usepackage{kotex}
\usepackage{lineno,hyperref}
\modulolinenumbers[5]
\usepackage{epstopdf}
\usepackage{pifont}
\usepackage{latexsym}
\usepackage{graphics}
\usepackage{graphicx}
\usepackage{float}
\usepackage[dvips]{xy}
\xyoption{all}
\usepackage{ifpdf}
\usepackage{multirow}
\usepackage{subfigure}
\usepackage{amssymb, amsthm, amsmath}
\usepackage{hhline}
\usepackage{centernot}
\usepackage{verbatim, comment, color}
\usepackage{enumerate}
\usepackage{setspace}
\usepackage[normalem]{ulem}

\usepackage{lipsum}

\newtheorem{thm}{Theorem}
\newtheorem{prop}{Proposition}
\newtheorem{lem}{Lemma}
\newtheorem*{pf}{Proof}

\newtheorem{mydef}{Definition}

\newtheorem{remark}{Remark}
\newtheorem{example}{Example}

\newcommand{\ba}{\begin{align}}
\newcommand{\ea}{\end{align}}  

\newcommand{\be}{\begin{equation}}
\newcommand{\ee}{\end{equation}}
\newcommand{\bea}{\begin{eqnarray}}
\newcommand{\eea}{\end{eqnarray}}
\newcommand{\barr}{\begin{array}}
\newcommand{\earr}{\end{array}}
\newcommand{\bn}{\begin{enumerate}}
\newcommand{\en}{\end{enumerate}}
\newcommand{\bi}{\begin{itemize}}
\newcommand{\ei}{\end{itemize}}
\newcommand{\bbbm}{\begin{pmatrix}}
\newcommand{\eeem}{\end{pmatrix}}

\newcommand{\p}{\partial}

\newcommand{\bbE}{{\bf E}}

\newcommand{\mathN}{{\mathbb N}}

\newcommand{\R}{{\mathbb R}}

\newcommand{\ep}{\epsilon}

\newcommand{\la}{\lambda}

\newcommand{\om}{\omega}

\newcommand{\ignore}[1]{}{}
\newcommand{\noin}{\noindent}

\newcommand{\nn}{\nonumber}

\newcommand{\q}{\quad}

\renewcommand{\subset}{\subseteq}

\renewcommand{\phi}{\varphi}
\renewcommand{\epsilon}{\varepsilon}
\newcommand{\cal}{\mathcal}

\DeclareMathOperator{\supp}{supp}

\numberwithin{equation}{section}
\numberwithin{theorem}{section}
\begin{document}
\title[uniqueness and characterization of local minimizers]{uniqueness and characterization of local minimizers for the interaction energy with mildly repulsive potentials}

\date{\today}



\author{Kyungkeun Kang, Hwa Kil Kim, Tongseok Lim, Geuntaek Seo}

\address{Kyungkeun Kang: Department of Mathematics \newline Yonsei University, Seoul 03722, Korea (Republic of)}
\email{kkang@yonsei.ac.kr}

\address{Hwa Kil Kim: Department of Mathematics Education \newline Hannam University, Daejeon 34430, Korea (Republic of)}
\email{hwakil@hnu.kr}

\address{Tongseok Lim: Institute of Mathematical Sciences \newline ShanghaiTech University, 393 Middle Huaxia Road, Pudong, Shanghai}
\email{TLIM@shanghaitech.edu.cn / TLIM0213@outlook.com}

\address{Geuntaek Seo: Department of Mathematics \newline Yonsei University, Seoul 03722, Korea (Republic of)}
\email{gtseo@yonsei.ac.kr}


\begin{abstract}
In this paper, we are concerned with local minimizers of an interaction energy governed by repulsive-attractive potentials of power-law type in one dimension.
We prove that sum of two Dirac masses is the unique local minimizer under the $\lambda-$Wasserstein metric topology with $1\le \lambda<\infty$,
provided masses and distance of Dirac deltas are equally half and one, respectively. In addition, in case of $\infty$-Wasserstein metric, we characterize
stability of steady-state solutions depending on powers of interaction potentials.
\end{abstract}

\maketitle


\section{Introduction and main results}

In this paper, we consider the following \textit{interaction energy}
\ba \label{energy}
 \rho \in \mathcal{P}(\R^d) \mapsto E(\rho) =\frac{1}{2} \iint V(x-y) d\rho(x)d\rho(y),
\end{align}
where  $\mathcal{P}(\R^d)$ denotes the space of Borel probability
measures on $\R^d$ and $V: \R^d \rightarrow (-\infty, \infty]$ is an
interaction potential. The interaction energy is in close connection
with the following continuity equation:
\begin{align}\label{pde}
\frac{\p \rho}{\p t} + \nabla \cdot (\rho \bold{v} )&=0,  \quad \bold{v}=-\nabla V \ast \rho, \q \text{where} \\
\nabla V \ast \rho(x)&=\int_{\R^{d}} \nabla V(x-y)d\rho(y). \nn
\end{align}
It is known  that the equation \eqref{pde} has the structure of a {\em gradient flow} of the interaction energy \eqref{energy} \cite{ags08, cdfls11, s17}.



The equation \eqref{pde}  appears in many applications including physics, biology, etc; see \cite{bcp97, blt06, bt04, chm14, cmv03, cmv06,  me99, t00} and the
references therein. There, $- \nabla V(x-y)$ is regarded as the force that a particle at $y$ exerts on a  particle at $x$.
Usually, the interaction potential $V$ depends only on the distance between particles.
Thus it is natural to consider radially symmetric interaction potentials of the form $V(x)=\om(|x|)$, where $\om$ is a function defined on $\R_+$.
Also, in most cases of applications, particles tend to repel each other in the short range and still want to remain cohesive as a whole.
Therefore $V$ is chosen  to be repulsive ($\om'(r)<0$) towards the origin and attractive ($\om'(r)>0$) towards infinity.
Thus it is typical to choose $\om$ to be decreasing on $(0,M)$ and increasing on $(M,\infty)$  with a unique minimum at $r=M$.

%

As natural candidates of steady states of \eqref{pde}, local and global minimizers of interaction energy \eqref{energy} and their uniqueness and structure, have received
substantial attention in recent years.
Regarding the various interaction potentials, 
there has been much interest in geometric properties of support of the minimizers  \cite{bclr13,  cfp17, ch17}.

For the case of purely attractive potentials, such as the Newtonian potential, the shape of the minimizer is well known as one point mass,
and  consequences of asymptotic dynamics have been widely studied \cite{bcl09, bl07, bll12, blr11}.
 On the other hand, in the case of repulsive-attractive potentials,  the geometry of 
 minimizers are sensitive to the precise form of the potential and they show various patterns \cite{bksu11, bkuv12, hkp13}. 
 In particular, local stability properties of steady states of \eqref{pde} for repulsive-attractive potentials have only been analyzed very recently \cite{bclr2013, fh13, fhk11}.

In this paper, we focus on the following interaction potential $V$,
\begin{equation}\label{pq}
V_{p,q}(x)=|x|^p /p - |x|^q / q, \q  p>q.
\end{equation}
This potential \eqref{pq}, called the power-law potential, is one of the most commonly considered potentials among repulsive-attractive potentials.
The term $-|x|^q / q$ is the repulsive one and $|x|^p/p$ is attractive one.
In this paper, our interest for $V_{p,q}$ is concentrated on the case $q\ge2$. In this case, $V_{p,q}$ is called {\em mildly repulsive}
\cite{bclr13, cfp17}.

There were early works towards understanding the steady states and asymptotic behavior of \eqref{pde} for the case of power-law potential \eqref{pq} (see $e.g.$ \cite{bclr2013, ch17, fh13, fhk11}).
In many cases, it is fairly general to assume spherical symmetry on a steady state solution, and then to find its exact form. The spherical symmetry assumption would often be supported by numerical simulations. 

However, in the case $q \ge 2$, it turns out the weakly repulsive force results in accumulation of particles \cite{bclr13,  cfp17}. This in particular implies that one cannot expect spherical symmetry of local minimizers. In such a situation where symmetry breaks, it is known to be hard to verify the uniqueness of minimizers or to precisely determine their structure.

 In this paper, we resolve these problems in one dimension. More precisely we show that for the problem \eqref{energy}
with $V_{p,q}$, not only  global minimizers but also local minimizers with respect to the Wasserstein metric $d_\la \  (1\leq \la <\infty)$ on $\cal P(\R)$  are uniquely characterized (up to translations) as the form $\frac{1}{2}\delta_{0} + \frac{1}{2} \delta_{1}$ if $p>>q \ge 2$ in Theorem  \ref{main}.

Furthermore, in Theorem \ref{sub}, we give a complete characterization of the measures of the form $\rho_m^{*}:=m \delta_0 + (1-m)\delta_1,~0<m<1$, with respect to the $d_\infty$-metric. That is, we determine  whether $\rho_m^{*}$ is a $d_\infty$-strict local minimizer or a saddle point of the energy \eqref{energy} for each and every $p>q\ge 2$. See also Remark 2 for an interpretation of the result in terms of the asymptotic stability of the evolution equation \eqref{pde}.

We note that the interaction energy is invariant under rigid motions, i.e. translation and rotation (which is the reflection on $\R$) of $\rho$. Hence, it is natural to regard two measures as the same if one is equal to the other via a rigid motion, and it is the viewpoint in this paper.  

 In order to precisely state our results about local energy minimizers, we need to discuss the topology on $\mathcal{P}(\R^d)$. For $1 \le \la \le \infty$ we will consider the $\la$- Wasserstein metric $d_\la$ on $\mathcal{P}(\R^d)$, which is defined as
\begin{align}
d_\la(\mu, \nu)&= \bigg{[}\inf_{\pi \in \Pi(\mu,\nu)}\iint |x-y|^\la d\pi(x,y) \bigg{]}^{1/\la} \q \text{for} \q 1 \le \la < \infty, \nn \\
d_\infty(\mu, \nu) &= \inf_{\pi \in \Pi(\mu,\nu)} \sup_{(x,y)\in \supp(\pi)} |x-y|, \nn
\end{align}
where $\mu, \nu \in \mathcal{P}(\R^d)$, and $\Pi(\mu,\nu) \subset \mathcal{P}(\R^d \times \R^d)$ consists of all couplings (also called transport plans) of $\mu, \nu$. The metric $d_\la$ is well-defined on $\mathcal{P}_\la(\R^d)$, the set of probability measures with finite moments of order $\la$. See \cite{ags08, bclr13, v03} for more details. Note that local minimizers in the $d_\la$- topology are automatically local minimizers in the $d_\infty$- topology.

Now we are ready to state our results. 
\begin{thm}\label{main} Let $d=1$, $q \ge 2$ and  $1 \le \la < \infty$. There exists $p^*=p^*(q)$ such that if $p > p^*$, then 
$\rho^*= \frac{1}{2}(\delta_0 + \delta_1)$ is the unique 
$d_\la$- local minimizer up to translation for the energy \eqref{energy}. In particular, $\rho^*$ is the unique global minimizer.
\end{thm}

\begin{remark}\label{remark1}
 It is clear  to see, by taking the infimum of all such $p^*$ in the theorem, that $p^*=p^*(q)$ is a unique function of $q\ge2$ which may be interpreted as a threshold for the phase transition of solutions of \eqref{energy}. A natural question is whether $p^*(q)=q$ or not. The following example indicates that it is not the case in general.
\end{remark}
\begin{example}\label{example}\footnote{We thank to Donghui Kim for providing the computation.} Let $p=2.5, q=2.1$, and let
$$\rho = (0.420137) \delta_0 + (0.159726) \delta_{0.548674} + (0.420137) \delta_{1.09735}.$$
Numeric computations confirm that $E(\rho) \approx - 0.0192448$, while $E(\rho^*) \approx -0.0190476$. This indicates that $\rho^*=\frac{1}{2}(\delta_0 + \delta_1)$ is not a global minimizer for some $p>q>2$. On the other hand if $q$ is large enough, for all $p>q$ simulations always seem to converge to $\rho^*$. This leads us to conjecture:
\end{example}
\noindent{\bf Conjecture.} There exists $q^* \ge 2$ so that for $q > q^*$, we have $p^*(q)=q$.

\vspace{2mm}
Related to this, see Section \ref{FO} for further remarks. 
\vspace{2mm}

A probability measure $\rho \in \cal P(\R^d)$ is called a {\em steady state} (see \cite{bclr2013}) if 
\[ -(\nabla V_{p,q} * \rho)(x)=0 \q \text{for all } \ x \in \supp(\rho).\]

\begin{mydef} Let $\rho \in \cal P(\R^d)$ be a steady state.\\
(1) $\rho$ is a $d_\infty$-strict local minimizer of \eqref{energy} if there exists $\ep >0$ such that for all $\rho' \in \cal P(\R^d)$ with $d_\infty(\rho, \rho') < \ep$ we have $E(\rho) \le E(\rho')$, and moreover $E(\rho) = E(\rho')$ if and only if  $\rho'$ is a translation of $\rho$. In other words, $\rho$ is the unique minimizer in its own small neighborhood.\\
(2) $\rho$ is a $d_\infty$-saddle point of \eqref{energy} if for every $\ep >0$ $\rho$ is neither a minimizer nor a maximizer in the $\ep$-ball $\{\mu \ | \ d_\infty(\mu, \rho) < \ep$\}.
\end{mydef}
The measures $m\delta_0 + (1-m)\delta_1$, $0<m<1$, have received particular interest, and one reason may be that they are steady states. In the following theorem, we characterize them with respect to the $d_\infty$- metric.

\begin{thm}\label{sub} Let $d=1$ and let $\rho_m^*=m\delta_0 + (1-m)\delta_1$, $0<m<1$.\\
(1) If $p>q> 2$, for every $m \in (0,1)$, $\rho_m^*$ is a $d_\infty$-strict local minimizer.\\
(2) If $p>3, q=2$, for every $m \in (\frac{1}{p-1}, \frac{p-2}{p-1})$, $\rho_m^*$ is a $d_\infty$-strict local minimizer.\\
(3) If $p>3, q=2$, for every $m \in (0, \frac{1}{p-1}] \cup [\frac{p-2}{p-1},1)$, $\rho_m^*$ is a $d_\infty$-saddle point.\\
(4) If $3>p>q= 2$, for every $m \in (0,1)$, $\rho_m^*$ is a $d_\infty$-saddle point.\\
(5) If $p=3, q=2$, for every $m \in (0, \frac{1}{2}) \cup (\frac{1}{2},1)$, $\rho_m^*$ is a $d_\infty$-saddle point.\\
(6) If $p=3, q=2$, $\rho_{\frac{1}{2}}^*$ is a $d_\infty$-strict local minimizer. 
\end{thm}
Upon the completion of this research, we realized some of the above results could follow by other works, e.g. it seems (1) and (2) could be a consequence of Theorem 3.1 in \cite{fr10} which pursues the dynamic point of view \eqref{pde}. Nevertheless, to the best of the authors' knowledge and search, we believe several other cases treated in this theorem are novel. In addition, our variational resolution may provide a unified viewpoint, and even for the cases (1), (2) our approach may give a new insight, as it clearly shows that certain  quadratic estimates can be applied.

But there are more subtle cases where such quadratic estimates are no longer available, and we need more careful investigation. We note there are two borderline cases, namely $p>3, q=2$ and $m=\frac{1}{p-1}$ in (3), and $p=3, q=2$ and $m=\frac{1}{2}$ in (6).  It is interesting to see that they exhibit the opposite characteristics. In particular, the resolution of (6) calls for the following estimate, which may be of independent interest.  
\begin{prop}\label{estimate} For $n  \in \mathbb{N}, \,  M >0$, assume $X, Y$ are i.i.d. random variables  with $|X| \le  (\frac{2n-1}{3M})^{\frac{1}{2n-1}}$ and  $\bbE[X^{2j-1}]=0$ for all $j=1,..., n$. Then
\begin{align} 
\bbE[|X-Y|^{2n+1}]-2\bbE[X^{2n+1}] -M  (\bbE[X^{2n}])^2  \ge 0. \nn
\end{align}
\end{prop}
We note that the proposition is sharp in the following sense: if $2\bbE[X^{2n+1}]$ is replaced with $C\bbE[X^{2n+1}]$ for any $C>2$, or $(\bbE[X^{2n}])^2$ with $\bbE[X^{4n}]$, the estimate may no longer hold. One can check this e.g. when $n=1$  by direct computation with centered random variables $X$ attaining only two real  values.
\begin{remark}[An interpretation of Theorem \ref{sub} in terms of stability]  It is known that a gradient flow, a solution $\{\rho_t\}_{t \ge 0}$ to the evolution equation \eqref{pde}, exhibits energy decay, i.e. $t \mapsto E(\rho_t)$ is a nonincreasing function of $t$. This implies that, under $d_\infty$-topology, every asymptotically stable state is a local minimizer. We therefore conclude that every saddle point in Theorem 2 is not asymptotically stable. This observation complements the results of \cite{fr10}.
\end{remark}
This paper is organized as follows. In section 2 we prove Theorem \ref{main}, and in section 3  we prove Theorem \ref{sub} along with Proposition \ref{estimate}. Lastly, in section 4 we give further remarks.

\section{Proof of Theorem \ref{main}}
Throughout the paper we fix a $q \ge 2$ and let $p >q$. For convenience, $V_{p,q}$ will be denoted by $V$. 
As $V$ is radial, by abusing notation we may regard $V$ as a function on $\R_+$, and define the following. Let $r>0$ be the unique inflection point of $V$ (i.e. $V''(r)=0$), and $R>0$ be the unique zero of $V$ (i.e. $V(R)=0$). Let $l=R-r$. It is easy to find
\ba \label{powerlaw}
r= \bigg{(}\frac{q-1}{p-1}\bigg{)}^{\frac{1}{p-q}}, \q R=\bigg{(}\frac{p}{q}\bigg{)}^{\frac{1}{p-q}}. \nn
\end{align}
As $p \to \infty$, we see that $r \nearrow 1$, $R\searrow 1$, and $l \searrow 0$. Keep in mind that $r, R, l$ are functions of $p,q$ (or functions of $p$, as we fixed $q$).

We start with the following lemma.
\begin{lem}\label{0th lemma}
Let $1 \le \la < \infty$, and let $\mu \in \mathcal{P}_\la(\R^d)$ be a $d_\la$- local minimizer with $E(\mu)<\infty$.  Then \emph{diam(supp($\mu$))$\leq R$}. 
\end{lem}
\begin{pf}
Choose any two points $z_1, z_2$ in $\supp(\mu)$. Let us define
$$
\nu_n:=\frac{1}{\mu(B(z_1,\frac{1}{n}))} \mu | _{B(z_1,\frac{1}{n})} - \frac{1}{\mu(B(z_2,\frac{1}{n}))} \mu | _{B(z_2,\frac{1}{n})}
$$
for all large $n \in \mathN$. Since $\mu$ is a local minimizer w.r.t $d_\la$-metric, there exists $\eta>0$ such that $d_\la(\mu, \rho)\leq  \eta$ implies $E(\mu)\leq E(\rho)$. It is clear that for each $n \in \mathN$, there exists $\gamma=\gamma(n)$ such that for all $0 <\ep < \gamma$, $\mu \pm \ep \nu_n \in \mathcal{P}_\la(\R^d)$ and $d_\la(\mu, \mu \pm \ep \nu_n) \leq \eta$. We observe that
\begin{align}
&E(\mu \pm \ep \nu_n) -  E(\mu) \nn \\
= \,\,&\pm \ep \iint V(x-y) d\mu(x) d\nu_n (y) + \frac{\ep^2}{2}\iint V(x-y) d\nu_n (x) d\nu_n (y) \nn \\
:=\,\,&\pm \ep I_1(n) + \ep^2 I_2(n) \geq 0. \nn
\end{align}
Note that $I_1, ~~I_2$ are finite. Now observe that since $\ep>0$ is arbitrary, we must have $I_1(n)=0.$ Then again by the above inequality, we obtain $I_2(n) \geq 0$. Letting $n \rightarrow \infty$, we have
$$
\lim_{n\rightarrow \infty} I_2(n)=-V(z_1 -z_2) \geq 0.
$$
This implies $V(z_1 - z_2) \leq 0$ for any $z_1, z_2$ in $\supp(\mu)$. Hence
by the property of $V$ we have $|z_1-z_2| \le R$, therefore ${\rm diam}(\supp(\mu)) \leq R$.
\qed
\end{pf}

From now on we will confine ourselves to the one-dimension $d=1$. By translation, we will always assume for any $d_\la$- local minimizer $\rho$,
\ba
\inf (\supp(\rho))=0, \q  \text{hence}  \q \supp(\rho) \subset [0,R]. \nn
\end{align}

\begin{lem}\label{first lemma}
If $l < r$, $\supp(\rho) \cap (l,r) = \emptyset$ for any $d_\la$- local minimizer
 $\rho$.
 \end{lem}
 \begin{pf} 
  For any $x \in (l,r)$ and $y \in [0,R]$ we have $|x-y| <r$, hence 
 \begin{equation}\label{ccc}
 \frac{d^2}{dx^2}V(x-y) <0.
 \end{equation} 
 This implies, if $x \in \supp(\rho) \cap (l,r)$ and $\rho$ is a local minimizer,  the energy must strictly decrease  if we slightly translate the mass of $\rho$ around $x$. More precisely, for small $\ep >0$, let $\rho\big|_{(x-\ep,x+\ep)}$ be the restriction of $\rho$ on $(x-\ep,x+\ep)$, and let $\eta_t$ be the translation of $\rho\big|_{(x-\ep,x+\ep)}$ by $t$. Then
\[
t \mapsto E(\rho - \rho\big|_{(x-\ep,x+\ep)} + \eta_t) \text{ is strictly concave around $t=0$ by \eqref{ccc}.} \]
But this contradicts to the assumption that  $\rho$ is a local minimizer.  \qed
\end{pf}

From now on assume $l <r$, which is the case if $p$ is large enough. Given a $d_\la$- local minimizer $\rho$ (for a given $p>q\ge2$), we will denote
\[
\rho_0 := \rho \big{|}_{[0,l]}, \q \rho_1 := \rho \big{|}_{[r,R]}, \q \text{so that} \q \rho= \rho_0 + \rho_1. \]
Denote $|\mu| := \mu(\R)$ for a positive measure $\mu$ on $\R$.

\begin{lem} \label{balance}  For any $d_\la$- local minimizer $\rho$, we have
\ba
\frac{{\rm max}(|\rho_0|, |\rho_1|)}{{\rm min}(|\rho_0|, |\rho_1|)} \le  \frac{V(1)}{V(1-3l) +\frac{1}{q}l^q},  \nn
\end{align}
whenever $V(1-3l) +\frac{1}{q}l^q<0$. Hence, $|\rho_0| \to 1/2$, $|\rho_1| \to 1/2$ as $p \to \infty$.
 \end{lem}
\begin{pf} By translation and reflection, suppose $|\rho_0| > |\rho_1|$ without loss of generality. Take a small $\ep \in (0,l)$, and recall $\inf (\supp(\rho))=0$. Let $\eta_0 :=  \rho_0\big|_{[0, \ep)}$, let $\eta_1$ be the translation of $\eta_0$ by $r-\ep$ (so that $\eta_1$ is concentrated on $[r-\ep,r)$), and  $\rho_0':= \rho_0 - \eta_0$. Then $\rho = \rho_0' + \eta_0+ \rho_1$. Define  $\rho^* := \rho_0' + \eta_1 + \rho_1$. Observe
  \ba
 E(\rho^*) - E(\rho) &= \iint V(x-y) d\eta_1(x)d\rho_0'(y) + \iint V(x-y) d\eta_1(x)d\rho_1(y)\nn \\
 & - \iint V(x-y) d\eta_0(x)d\rho_0'(y) - \iint V(x-y) d\eta_0(x)d\rho_1(y). \nn
\end{align}
We estimate as follows; recall $|\eta_0|=|\eta_1|$. 
  \ba
& \iint V(x-y) d\eta_1(x)d\rho_0'(y)  - \iint V(x-y) d\eta_0(x)d\rho_0'(y)   \nn \\
&\le  |\eta_1|\int V(r-\ep-y)  d\rho_0'(y)  - |\eta_0| \int V(0-y)  d\rho_0'(y) \nn \\
& \le |\eta_0|  \int [V(r-\ep-y) +\frac{1}{q}|y|^q] d\rho_0'(y) \q (\text{since } -V(x) \le \frac{1}{q}|x|^q) \nn \\
& \le  |\eta_0| |\rho_0'| ( V(r-\ep-l) +\frac{1}{q}l^q) \q  (\text{since } \ 0 \le y \le l) \nn \\
& \le |\eta_0| (|\rho_0|-|\eta_0|) ( V(1-3l) +\frac{1}{q}l^q) \q  (\text{since } \ r-\ep-l \ge 1-3l). \nn
\end{align}
 Next, we estimate
   \ba
& \iint V(x-y) d\eta_1(x)d\rho_1(y)  - \iint V(x-y) d\eta_0(x)d\rho_1(y)   \nn \\
&\le |\eta_1| \int [V(r-y) -  V(0-y)] d\rho_1(y)  \nn \\
&\le |\eta_1| \int [-  V(y)] d\rho_1(y) \q (\text{since } |r-y| \le l \text{ implies } V(r-y) \le 0) \nn \\
& \le -|\eta_1| |\rho_1| V(1). \q (\text{since } V(1) \le V(x) \q \text{for all } \, x) \nn
\end{align}
Combining, we get
\ba
&E(\rho^*) - E(\rho)\nn\\
& \le |\eta_0| \bigg{(}|\rho_0| ( V(1-3l) +\frac{1}{q}l^q) -  |\rho_1| V(1)\bigg{)} - |\eta_0|^2 \bigg{(} V(1-3l) +\frac{1}{q}l^q\bigg{)}. \nn
\end{align}
Now if $\rho(\{0\})=0$, $\ep \to 0$ yields $|\eta_0| \to 0$. This implies that\\  $|\rho_0| ( V(1-3l) +\frac{1}{q}l^q) -  |\rho_1| V(1) \ge 0$, otherwise we get $E(\rho^*) - E(\rho)<0$ for small $\ep$, a contradiction to the local minimality of $\rho$. We conclude 
\begin{align}\label{ratio}
\frac{|\rho_0|}{|\rho_1|} \le \frac{V(1)}{V(1-3l) +\frac{1}{q}l^q}.  
\end{align}
If $\rho(\{0\})>0$, $\ep \to 0$ may not yield $|\eta_0| \to 0$. But in this case we take $\ep  < \rho_0(\{0\})$, and let $\eta_0 = \ep \delta_0$, $\eta_1 = \ep \delta_1$, $\rho_0'= \rho_0 - \eta_0$, and  define $\rho^*$ as before. By following the similar estimates and taking $\ep \to 0$ we again obtain \eqref{ratio} (In fact, we can get a slightly stronger inequality $\frac{|\rho_0|}{|\rho_1|} \le \frac{V(1)}{V(1-l) +\frac{1}{q}l^q}$.)

 In particular, we see that $p \to \infty$ yields $l \to0$ and hence $\frac{|\rho_0|}{|\rho_1|} \to 1$.
 \qed
 \end{pf}
\begin{lem}\label{localbound1}
Let $q \ge 2$, $k>0$. If  $p\geq q+k+1$, then $(x-1)V'(x) \ge k(x-1)^2$ for all $x \in [c, \infty)$, where $c \in (0,1)$ is the solution to the equation $x^{q+k}-x^{q-1}=k(x-1)$.
 \end{lem}
\begin{pf}
Let $g(x)=V'(x)=x^{p-1}-x^{q-1}$. We have
\ba
g'(x)&= (p-1)x^{p-2}- (q-1)x^{q-2}, \nn \\
g''(x)&=(p-1)(p-2)x^{p-3}- (q-1)(q-2)x^{q-3}. \nn
\end{align}
We see that $g$ has decreasing-increasing \& concave-convex shape on $\R_+$ and has a unique inflection point in $(0,1)$. In particular, $g$ is convex on $[1,\infty)$. Since $g(1)=0$ and $g'(1)=p-q$, for all $p \ge q+k$ we have $g(x) \ge k(x-1)$ on $[1,\infty)$, hence $(x-1)V'(x) \ge k(x-1)^2$ on $[1,\infty)$.\\
Next, we compare $g(x)$ and $k(x-1)$ when $0 \le x\le 1$. By the shape of $g$ it is clear that if $g'(1) > k$, there exists a unique solution $x_{p,k}$ in $(0,1)$ to the equation $g(x)=k(x-1)$, and $g(x) \le k(x-1)$ on $[x_{p,k}, 1]$. Since $g(x) \leq x^{q+k}-x^{q-1}$ for $p\geq q+k+1$, we have $x_{p,k} \le c$. Hence $g(x)\leq k(x-1)$ for all $x \in [c, 1]$. By multiplying $(x-1)$ on both sides, the lemma follows.  \qed
\end{pf}
Now we prove Theorem \ref{main}.
\begin{pf}[Theorem 1]
Let $\rho= \rho_0 +  \rho_1$ be a $d_\la$- local minimizer. Let $m_0=|\rho_0|$, $m_1=|\rho_1|$. We will consider the linear contraction of $\rho_0$ to $m_0\delta_0$ and $\rho_1$ to $m_1 \delta_1$ respectively. That is, we linearly transfer the mass $\rho_0$ on  $[0,l]$ to $0$ and $\rho_1$ on $[r,R]$ to $1$ as $t$ goes  from $0$ to $1$. Let $\rho_t$ be the contraction of $\rho$ at time $t$, so that $\rho_0=\rho$, $\rho_1=m_0\delta_0 + m_1 \delta_1$. Then
\ba
E(\rho_t) &= \iint V \big{(} t(1-y+x) +(y-x)\big{)} d\rho_0(x)d\rho_1(y) \nn \\
& + \sum_{i=0}^1\frac{1}{2} \iint  V\big{(}(1-t)(y-x)\big{)} d\rho_i(x)d\rho_i(y).\nn
\end{align}
$\frac{d}{dt}E(\rho_t)\big{|}_{t=0} \ge 0$ as $\rho$ is a local minimizer. We differentiate and obtain
\ba
\frac{d}{dt}E(\rho_t)\bigg{|}_{t=0} &= \iint (1-y+x)V'(y-x) d\rho_0(x)d\rho_1(y) \nn \\
& - \sum_{i=0}^1\frac{1}{2} \iint  (y-x)V'(y-x) d\rho_i(x)d\rho_i(y).\nn
\end{align}
Note that $xV'(x) = x^p-x^q \ge -x^q \ge -x^2$ on $[0,1]$ as $q \ge 2$. In Lemma \ref{localbound1}, take $k >1$ and $p^*$ sufficiently large such that $y-x > c$ for all $p > p^*$ and $x \in \supp(\rho_0), y \in \supp(\rho_1)$. Then the lemma implies
\ba
&\frac{d}{dt}E(\rho_t)\bigg{|}_{t=0}  \nn \\
&\le -k \iint (y-x-1)^2 d\rho_0(x)d\rho_1(y)+ \sum_{i=0}^1\frac{1}{2} \iint  (y-x)^2 d\rho_i(x)d\rho_i(y).\nn
\end{align}
Now we compute the integrals. To see the computation more clearly, we adapt the following probabilistic notation: let $X$, $Y$ be random variables whose laws (distributions) are the probability measures $m_0^{-1}\rho_0$, $m_1^{-1}\rho_1$ respectively. Recall $\supp(\rho_0) \subset [0,l]$, $\supp(\rho_1) \subset [r,R]$. We compute
\ba
&-k \iint (y-x-1)^2 d\rho_0(x)d\rho_1(y) \nn \\
& -k m_0m_1\iint (y-x-1)^2 \, d[m_0^{-1}\rho_0](x)d[m_1^{-1}\rho_1](y) \nn \\
 & = -km_0m_1 [\bbE(X^2) +\bbE(Y^2) +1 + 2\bbE(X)-2\bbE(Y) -2\bbE(X)\bbE(Y)] \nn \\
 & = -km_0m_1 [{\rm Var}(X) +{\rm Var}(Y) \nn \\
& \q\q+\bbE(X)^2 + \bbE(Y)^2 +1 + 2\bbE(X)-2\bbE(Y) -2\bbE(X)\bbE(Y)] \nn \\
 & = -km_0m_1 [ {\rm Var}(X) +{\rm Var}(Y) + (\bbE(Y) -\bbE(X)-1)^2], \nn
\end{align}
where ${\rm Var}(X) = \bbE(X^2) - \bbE(X)^2$ is the variance of $X$. Next, we compute
 \ba
 &\frac{1}{2} \iint  (y-x)^2 d\rho_0(x)d\rho_0(y) \nn \\
& =\frac{m_0^2}{2} \iint  (y-x)^2 \,d[m_0^{-1}\rho_0](x)d[m_0^{-1}\rho_0](y) \nn\\
& = \frac{m_0^2}{2}(2\bbE(X^2) -2\bbE(X)^2) = m_0^2{\rm Var}(X). \nn
 \end{align}
Similarly,
$$\frac{1}{2} \iint  (y-x)^2 d\rho_1(x)d\rho_1(y) = m_1^2{\rm Var}(Y).$$
Hence we get
\ba
\frac{d}{dt}E(\rho_t)\bigg{|}_{t=0} \le & (-km_0m_1 +m_0^2){\rm Var}(X) + (-km_0m_1 +m_1^2){\rm Var}(Y) \nn \\
& -km_0m_1 (\bbE(Y) -\bbE(X)-1)^2. \nn
\end{align}
Since $m_0, m_1 \to 1/2$ as $p \to \infty$ by Lemma \ref{balance}, $k >1$ implies that there exists $p^*$ such that for all $p > p^*$, we have
$${\rm Var}(X) = 0, \q {\rm Var}(Y)=0, \q \bbE(Y) -\bbE(X) =1,$$
since otherwise we have $\frac{d}{dt}E(\rho_t)\bigg{|}_{t=0}<0$, a contradiction to the fact that $\rho$ is a local minimizer. Hence $\rho$ is of the form $\rho= m\delta_0 + (1-m)\delta_1$. In this case $E(\rho) = m(1-m)V(1)$, which is minimized when $m=1/2$. Again, the fact that $\rho$ is a $d_\la$- local minimizer implies that $m=1/2$. We conclude that for all $p > p^*$, $\frac{1}{2}(\delta_0 +\delta_1)$ is a unique $d_\la$- local (hence global) minimizer for the energy \eqref{energy}. \qed
\end{pf}

\section{Proof of Theorem \ref{sub}}
Firstly, we prove Proposition \ref{estimate}.
\begin{pf}[Proposition \ref{estimate}]
Let $\mu \in \mathcal{P}(\R)$ be the distribution of $X$ (and so $Y$) which is concentrated on $[-c,c]$ for some $c>0$. Observe
\begin{align} &\bbE[|X-Y|^{2n+1}]-2\bbE[X^{2n+1}] -M  (\bbE[X^{2n}])^2 = \bbE f(X,Y), \text{ where}  \nn \\
&f(X,Y):=|X-Y|^{2n+1}-X^{2n+1}-Y^{2n+1}-MX^{2n}Y^{2n}. \nn
\end{align}
We may assume that $\mu$ is absolutely continuous with respect to Lebesgue measure, as the result for general $\mu$ can be obtained by approximation by absolutely continuous measures.
By symmetry of $f$, we can rewrite
\begin{align} \bbE f(X,Y) &= \iint f(x,y)d\mu(y)d\mu(x) =2 \int_{-c}^{c} \int_{-c}^x g(x,y)d\mu(y)d\mu(x),    \nn
\end{align}
where 
\begin{align}
g(x,y) &= (x-y)^{2n+1}-x^{2n+1}-y^{2n+1}-M x^{2n}y^{2n} \nn \\
&=-2y^{2n+1}+\sum_{k=1}^{2n}(-1)^{2n+1-k}\binom{2n+1}{k}x^k y^{2n+1-k}-M x^{2n}y^{2n}. \nn
\end{align}

Let $\nu:= \mu \otimes \mu$ be the tensor product, and define the following regions
\begin{align}
A&=\{(x,y) \ | \ 0 \le x \le c,\  0 \le y \le x\}, \nn \\
A' &= \{(x,y) \ | \ 0 \le x \le c,\  x \le y \le c\},   \nn \\
B&= \{(x,y) \ |  -c \le x \le 0, \ -c \le y \le x\}, \nn \\
B'&= \{(x,y) \ |  -c \le x  \le 0, \ x \le y \le 0\}, \nn \\
C&=\{(x,y) \ | \ 0 < x \le c, \ -c \le y < 0\}.  \nn
\end{align}
Then we can decompose
\[
\frac{1}{2}\bbE f(X,Y) = \iint_{A \cup B\cup C} g \, d\nu
\]
\[
=-\iint_{B\cup C} 2y^{2n+1}\, d\nu + \iint_A (2xy^{2n}-2y^{2n+1}) \,d\nu + \iint_A (2n-1)xy^{2n} \,d\nu
\]
\[
+ \iint_{B \cup C} (2n+1)xy^{2n}\, d\nu -\iint_{A \cup C} (2n+1)x^{2n}y\,d\nu - \iint_B (2n+1)x^{2n}y\,d\nu
\]
\[
+\sum_{k=2}^{2n-1}(-1)^{2n+1-k}\binom{2n+1}{k}\iint_{A \cup B \cup C}x^k y^{2n+1-k}d\nu- \iint_{A \cup B \cup C} Mx^{2n}y^{2n}\,d\nu.
\]
Notice $-2y^{2n+1} \ge 0$ on $B\cup C$ and $2xy^{2n}-2y^{2n+1} \ge 0$ on $A$. Also, see
\[ \iint_{B \cup C} xy^{2n}\, d\nu \ge \iint_{x \in \R, y < 0} xy^{2n}\, d\nu = 0,
\]
since $\bbE[X]=0$. Similarly $\iint_{A \cup C}x^{2n}y\,d\nu \le 0$. We further claim:
\begin{equation}\label{KKLS-20}
\sum_{k=2}^{2n-1}(-1)^{2n+1-k}\binom{2n+1}{k}\iint_{A \cup B \cup C}x^k y^{2n+1-k}d\nu\ge 0.
\end{equation}
Indeed,
\begin{align}
&\sum_{k=2}^{2n-1}(-1)^{2n+1-k}\binom{2n+1}{k}\iint_{A \cup B \cup C}x^k y^{2n+1-k}d\nu \nn \\
&=\sum_{k=2}^{n}(-1)^{2n+1-k}\binom{2n+1}{k}\iint_{A \cup B \cup C}x^k y^{2n+1-k}d\nu
\nn \\
&+\sum_{k=n+1}^{2n-1}(-1)^{2n+1-k}\binom{2n+1}{k}\iint_{A \cup B \cup C}x^k y^{2n+1-k}d\nu \nn \\
&=\sum_{k=2}^{n}(-1)^{k+1}\binom{2n+1}{k}\iint_{A \cup B \cup C} (x^k y^{2n+1-k}- x^{2n+1-k}y^k) d\nu, \nn
\end{align}
where we used the change of variables. 
Suppose $k$ is even. We then note
\[
\iint_{A \cup B \cup C} (x^k y^{2n+1-k}- x^{2n+1-k}y^k) d\nu\le 0.
\]
Indeed, since $2n+1-k$ is odd,
\begin{align}
&\iint_{A \cup B \cup C} x^k y^{2n+1-k}d\nu =\iint_{A \cup C} x^k y^{2n+1-k}d\nu+\iint_{B} x^k y^{2n+1-k}d\nu \nn \\
&\le \iint_{x>0, y \in \R} x^k y^{2n+1-k}d\nu+\iint_{B} x^k y^{2n+1-k}d\nu = \iint_{B} x^k y^{2n+1-k}d\nu \le 0, \nn
\end{align}
where $\bbE[X^{2n+1-k}]=0$ is used.
Moreover,
\begin{align}
&\iint_{A \cup B \cup C} x^{2n+1-k}y^k d\nu=\iint_{B \cup C} x^{2n+1-k}y^k d\nu+\iint_{A} x^{2n+1-k}y^k d\nu \nn \\
&\ge\iint_{x \in \R, y <0} x^{2n+1-k}y^k d\nu+\iint_{A} x^{2n+1-k}y^k d\nu = \iint_{A} x^{2n+1-k}y^k d\nu\ge 0. \nn
\end{align}
On the other hand if $k$ is odd, then $2n+1-k$ is even, and we can similarly obtain
\[
\iint_{A \cup B \cup C} (x^k y^{2n+1-k}- x^{2n+1-k}y^k) d\nu\ge 0,
\]
where we used $\bbE[X^{k}]=0$. This implies \eqref{KKLS-20}, and hence we deduce
\begin{align}\label{ineq}
&\frac{1}{2}\bbE f(X,Y)  \nn \\
& \ge  \iint_A (2n-1)xy^{2n} \,d\nu - \iint_B (2n-1)x^{2n}y\,d\nu - \iint_{A \cup B \cup C} Mx^{2n}y^{2n}\,d\nu. \nn
\end{align}
We then use the symmetry with respect to the diagonal to deduce 
\begin{align}
\iint_A xy^{2n} \,d\nu &= \frac{1}{2} \bigg( \iint_A xy^{2n} \,d\nu  + \iint_{A'} x^{2n}y \,d\nu \bigg) \nn \\
&\ge \frac{1}{2c^{2n-1}} \iint_{A \cup A'} x^{2n}y^{2n} \, d\nu =  \frac{1}{2c^{2n-1}} \bigg(\int_{\R_+} x^{2n} d\mu(x)\bigg)^2, \nn
\end{align}
where we used $|X|\le c$, $|Y| \le c$. Similarly,
\[
- \iint_B x^{2n}y\,d\nu \ge  \frac{1}{2c^{2n-1}} \bigg(\int_{\R_-} x^{2n} d\mu(x)\bigg)^2.
\]
On the other hand, we observe
\[
\iint_{A \cup B \cup C} x^{2n}y^{2n} d\nu \le \iint_{A \cup A'}x^{2n}y^{2n} d\nu + \iint_{B \cup B'}x^{2n}y^{2n} d\nu + \iint_{C} x^{2n}y^{2n} d\nu,
\]
and moreover
\[
\iint_{C} x^{2n}y^{2n} d\nu \le  \frac{1}{2}\bigg(\int_{\R_+} x^{2n} d\mu(x)\bigg)^2 +  \frac{1}{2}\bigg(\int_{\R_-} x^{2n} d\mu(x)\bigg)^2.
\]
Thus shows that $\bbE f(X,Y) \ge 0$  if we set $c \le (\frac{2n-1}{3M})^{\frac{1}{2n-1}}$. \qed
\end{pf}
Now we prove Theorem \ref{sub}.

\begin{pf}[Theorem \ref{sub}] Let $\rho^*=m\delta_0 + (1-m)\delta_1$ for some $0<m<1$. Let $m_0=m, m_1=1-m$. First of all, $\rho^*$ is certainly not a local maximizer, since $E(\rho^*) < E(m\delta_x + (1-m)\delta_1)$ for any $x \neq 0$, $x \neq 2$. \\

\noin$\bullet$ \emph{(Case} $p>q>2$\emph{)} Since $V''(1)=p-q$ and $V''(0)=0$, we note that
\begin{equation}\label{v''}
V''(1) > -V''(0).
\end{equation}
Let $a = \ep -V''(0)/2$, $b= V''(1)/2 - \ep$ for small $\ep >0$ so that $b > a > 0$. Then there exists $r_0=r_0(\ep)>0$ such that
$$V(x) \ge -ax^2 \text{ on } [-2r_0, 2r_0],$$
$$ V(x) -V(1) \ge b(x-1)^2 \text{ on } [1-2r_0, 1+2r_0].$$
We denote by $B_r(\rho^*)$ the ball of (small) radius $r$ and center $\rho^*$ in the $d_\infty$- metric. Note that any $\rho \in B_r(\rho^*)$ can be written as $\rho = \rho_0+\rho_1$, where
$$|\rho_0|=m_0, \q |\rho_1|=m_1, \q \supp(\rho_0) \subset [-r,r],\q  \supp(\rho_1) \subset [1-r,1+r].$$
We will directly compare the energies of $\rho$ and $\rho^*$. Observe

\ba
&E(\rho^*) - E(\rho)\nn \\
= &\iint [V(1)-V(x-y)] d\rho_0(x)d\rho_1(y)  - \sum_{i=0}^1\frac{1}{2} \iint V(x-y) d\rho_i(x)d\rho_i(y).\nn
\end{align}
Note that $y-x \in [1-2r, 1+2r]$ for $x \in [-r,r]$, $y \in [1-r,1+r]$, and $|y-x| \le 2r$ for $x,y \in [-r,r]$.
By employing probability notation as before, we compute
\ba
&\iint [V(1)-V(x-y)] d\rho_0(x)d\rho_1(y)  \nn \\
&= m_0m_1\iint [V(1)-V(y-x)] \, d[m_0^{-1}\rho_0](x)d[m_1^{-1}\rho_1](y) \nn \\
& \le -bm_0m_1 \iint (y-x-1)^2 \, d[m_0^{-1}\rho_0](x)d[m_1^{-1}\rho_1](y) \nn \\
 & = -bm_0m_1 [\bbE(X^2) +\bbE(Y^2) +1 +2\bbE(X)-2\bbE(Y) -2\bbE(X)\bbE(Y)] \nn \\
  & = -bm_0m_1 [ {\rm Var}(X) +{\rm Var}(Y) + (\bbE(Y) - \bbE(X) -1)^2] \nn
 \end{align}
Next, we compute
 \ba
& -\frac{1}{2} \iint V(x-y) d\rho_0(x)d\rho_0(y) \nn\\
& \le   \frac{am_0^2}{2} \iint (x-y)^2  \,d[m_0^{-1}\rho_0](x)d[m_0^{-1}\rho_0](y) \nn \\
& = \frac{am_0^2}{2}(2\bbE(X^2) -2\bbE(X)^2) = am_0^2{\rm Var}(X). \nn
 \end{align}
Similarly,
$$ -\frac{1}{2} \iint V(x-y) d\rho_1(x)d\rho_1(y) \le am_1^2{\rm Var}(Y).$$
 Combining the estimates, we get
\ba
 & E(\rho^*) - E(\rho)  \nn\\
 \le \, &(am_0^2-bm_0m_1){\rm Var}(X)+ (am_1^2-bm_0m_1){\rm Var}(Y) \nn\\
  & -  bm_0m_1(\bbE(Y) - \bbE(X) -1)^2. \nn
   \end{align}

Now if $p>q>2$, $V''(0)=0$ so $a \searrow 0$ as $\ep \searrow 0$. Hence for any given $m >0$, we get $am_0^2-bm_0m_1 <0$ and $am_1^2-bm_0m_1 <0$ for sufficiently small $\ep$. This implies $E(\rho^*) - E(\rho) <0$ unless ${\rm Var}(X)=0, {\rm Var}(Y)=0$ and $\bbE(Y) - \bbE(X) =1$. This proves that $\rho^*$ is a unique minimizer in $B_r(\rho^*)$ for any  $r$ with $0<r<r_0$.  \\

\noin$\bullet$ \emph{(Case} $p>3, q=2$, and $m \in (\frac{1}{p-1}, \frac{p-2}{p-1})$\emph)
We recall that \eqref{v''} is also valid since $V''(0)=-1$, $V''(1)=p-2$. Reminding that $a=\ep+1/2$ and $b=(p-2)/2 -\ep$, the inequalities $am_0^2-bm_0m_1 <0$ and $am_1^2-bm_0m_1 <0$ are equivalent to
$$\frac{1+2\ep}{p-1}<m_i <\frac{p-2 -2\ep}{p-1}, \q i=0,1.$$
As $\ep>0$ is arbitrary, this implies that for any given $m \in (\frac{1}{p-1}, \frac{p-2}{p-1})$, $\rho^*$ is a unique minimizer in $B_r(\rho^*)$ for any  $r$ with $0<r<r_0$.\\

\noin$\bullet$ \emph{(Case} $p>3, q=2$, and $m \in (0,\frac{1}{p-1}] \cup [\frac{p-2}{p-1},1 )$\emph{)}\\
Define  $\rho_x = \frac{m}{2} \delta_{-x} + \frac{m}{2} \delta_{x} + (1-m) \delta_1$, and consider
\begin{equation}\label{new energy}
E(x) := E(\rho_x) = \frac{m^2}{4}V(2x) + \frac{m(1-m)}{2} (V(1+x) +V(1-x)). \end{equation}
We find
\begin{equation}\label{E''}
E'(0)=0, \q E''(0) = -(p-1)m(m - \frac{p-2}{p-1})<0,
\end{equation}
provided that $m \in ( \frac{p-2}{p-1},1)$.
This implies that $\rho^*$ is a saddle point in each case. For the case that $m\in(0,\frac{1}{p-1})$, we take
$\rho_x=m\delta_0 +\frac{1-m}{2} \delta_{1-x} +\frac{1-m}{2} \delta_{1+x}$. Following the similar computations, we see that
\begin{equation}\label{another E''}
E''(0)=(1-m) \Big( m(p-2)-(1-m) \Big)<0,
\end{equation}
which is equivalent to the case $m \in (0,\frac{1}{p-1})$.\\

Now let us analyze the borderline case, that is,  $m= \frac{p-2}{p-1}$. The idea is to look at the third order Taylor expansion of $V$ at $x=1$. Observe that $V'(1)=0, V''(1)=p-2, V'''(1)=(p-1)(p-2)$ implies that given $0< b < \frac{(p-1)(p-2)}{6}$, there exists a small $c=c(b) >0$ such that
\begin{align}\label{third}
V(x)-V(1) \le \frac{p-2}{2}(x-1)^2 + b (x-1)^3 \q \text{on} \q x \in [1-c,1].
\end{align}
Using this, we make a similar estimate as before but in the opposite direction. Our aim is to find $\rho$ satisfying $E(\rho^*) > E(\rho)$. To this end, suppose $\rho$ is concentrated on $[0,c/2] \cup [1-c/2,1]$. By \eqref{third}, we have
\begin{align}
E&(\rho^*) - E(\rho)\nn \\
= &\iint [V(1)-V(y-x)] d\rho_0(x)d\rho_1(y)  - \sum_{i=0}^1\frac{1}{2} \iint V(y-x) d\rho_i(x)d\rho_i(y)
\nn \\
\ge & \iint [    -\frac{p-2}{2}(1-y+x)^2 +b(1-y+x)^3] d\rho_0(x)d\rho_1(y) \nn \\
& + \sum_{i=0}^1\frac{1}{4} \iint (y-x)^2 d\rho_i(x)d\rho_i(y) -  \sum_{i=0}^1\frac{1}{2p} \iint |y-x|^p d\rho_i(x)d\rho_i(y). \nn
\end{align}
As before, the three quadratic integrals sum up to
\ba
& \iint  -\frac{p-2}{2}(1-y+x)^2 d\rho_0(x)d\rho_1(y) + \sum_{i=0}^1\frac{1}{4} \iint (y-x)^2 d\rho_i(x)d\rho_i(y) \nn \\
& = (\frac{1}{2}m_0^2-\frac{p-2}{2}m_0m_1){\rm Var}(X)+ (\frac{1}{2}m_1^2-\frac{p-2}{2}m_0m_1){\rm Var}(Y) \nn\\
  & -  \frac{p-2}{2}m_0m_1(\bbE(Y) - \bbE(X) -1)^2,\nn
\end{align}
where $m_0=\rho([0,c/2])$, $m_1=\rho([1-c/2,1])$, $X \sim \frac{1}{m_0}\rho\big{|}_{[0,c/2]}$, $Y \sim \frac{1}{m_1}\rho\big{|}_{[1-c/2,1]}$. Notice that $m_0=\frac{p-2}{p-1}$ implies $\frac{1}{2}m_0^2-\frac{p-2}{2}m_0m_1=0$. Hence
\ba
&E(\rho^*) - E(\rho)\nn \\
& \ge  (\frac{1}{2}m_1^2-\frac{p-2}{2}m_0m_1){\rm Var}(Y)
  -  \frac{p-2}{2}m_0m_1(\bbE(Y) - \bbE(X) -1)^2   \nn \\
  &+ \iint b(1-y+x)^3 d\rho_0(x)d\rho_1(y)
   -  \sum_{i=0}^1\frac{1}{2p} \iint |y-x|^p d\rho_i(x)d\rho_i(y).   \nn
\end{align}
Since $\frac{1}{2}m_1^2-\frac{p-2}{2}m_0m_1 <0$, to achieve $E(\rho^*) - E(\rho) >0$ it is desired to set
${\rm Var}(Y)  =0$. This suggests that we may choose
$$\rho = (m-\ep)\delta_0 + \ep\delta_\eta + (1-m)\delta_1,$$
where  $\eta, \ep >0$ shall be chosen later. With this choice, note that $Y=1$, $\bbE(X)=\frac{\ep\eta}{m}$. Observing that
$$
\iint b(1-y+x)^3 d\rho_0(x)d\rho_1(y) =b(1-m)\ep\eta^3,
$$

$$
\sum_{i=0}^1\frac{1}{2p} \iint |y-x|^p d\rho_i(x)d\rho_i(y)=\frac{1}{p} (m-\ep) \ep\eta^{p} ,
$$
we obtain
\ba
&E(\rho^*) - E(\rho)\nn \\
& \ge
  -   \frac{1}{2}\bigg{(}\frac{p-2}{p-1}\bigg{)}^2\bigg{(}\frac{\ep\eta}{m}\bigg{)}^2
  +b(1-m)\ep\eta^3  -\frac{1}{p} (m-\ep) \ep\eta^{p} \nn \\
  & =   -   \frac{1}{2}\ep^2\eta^2
  +\frac{b}{p-1}\ep\eta^3 -\frac{p-2}{p(p-1)}\ep\eta^{p} +\frac{1}{p} \ep^2\eta^p \nn\\
  &= -   \frac{1}{2}\eta^{2\alpha+2}  + \frac{b}{p-1}\eta^{\alpha+3} -\frac{p-2}{p(p-1)}\eta^{\alpha+p} +\frac{1}{p} \eta^{2\alpha + p} \q \text{if} \q \ep=\eta^\alpha.\nn
\end{align}
By choosing $\alpha >1$ and $\eta$ sufficiently small, we get $E(\rho^*) - E(\rho)>0$, which implies that $\rho^*= \frac{p-2}{p-1}\delta_0 + \frac{1}{p-1}\delta_1$ is a saddle point. Notice  the case $m=\frac{1}{p-1}$ immediately follows by the reflection symmetry of the energy.\\

\noin$\bullet$ \emph{(Case} $3>p>q=2$\emph)
With $m\geq \frac{1}{2}$,  we can follow the same line of reasoning \eqref{new energy}, \eqref{E''}, \eqref{another E''} to get the desired result, and we omit the detail.  The case  $m \leq \frac{1}{2}$ is obviously obtained by reflection symmetry.\\

\noin$\bullet$ \emph{(Case} $p=3, q=2$, and $m \in (0,\frac{1}{2})\cup (\frac{1}{2}, 1)     $\emph) 
Again we can follow the same line of reasoning \eqref{new energy}, \eqref{E''}, \eqref{another E''}  and obtain the desired result.\\


\noin$\bullet$ \emph{(Case} $p=3, q=2$, and $m=\frac{1}{2}$\emph) Finally, we analyze the remaining borderline case $p=3, q=2, m=1/2$. In this case, since
$$V(x)-V(1) = \frac{1}{2}(x-1)^2 + \frac{1}{3} (x-1)^3 \q \text{on} \q x \ge 0$$
 the above estimates become exact, and we have
\ba
&E(\rho^*) - E(\rho)\nn \\
& =
  -  \frac{1}{8}(\bbE(Y) - \bbE(X) -1)^2   \nn \\
  &+ \iint \frac{1}{3}(1-y+x)^3 d\rho_0(x)d\rho_1(y)
   -  \sum_{i=0}^1\frac{1}{6} \iint |y-x|^3 d\rho_i(x)d\rho_i(y)   \nn
\end{align}
for all $\rho=\rho_0+\rho_1$ in some small $d_\infty$- neighborhood of $\rho^*=\frac{1}{2}\delta_0 + \frac{1}{2}\delta_1$. Recall that $|\rho_0|=|\rho_1|=1/2$, $\rho_0$ is concentrated in a neighborhood of $0$, say in $(-\ep, \ep)$ for some small $\ep >0$, $\rho_1$ is concentrated in $(1-\ep, 1+\ep)$, and $X,Y$ are independent random variables having distributions $ 2\rho_0,  2\rho_1$ respectively.  Let $u=\bbE(Y)-1$, $Z=Y-\bbE(Y)$. Let $X'$ be an independently and identically distributed random variable as $X$, and $Z'$ be i.i.d. as $Z$.

By translation invariance of the energy we can assume $\bbE(X)=0$, and this implies that $u$ is near $0$. With these notations we may rewrite
\ba
&24(E(\rho^*) - E(\rho))\nn \\
 =   &
  -3u^2
 -2\bbE[(u+Z-X)^3]
   -  \bbE[|X-X'|^3]  -  \bbE[|Z-Z'|^3]     \nn \\
    = &-3u^2 -2u^3-6u\bbE[(Z-X)^2]  \nn \\
   & -2\bbE[(Z-X)^3] -  \bbE[|X-X'|^3]  -  \bbE[|Z-Z'|^3].    \nn
\end{align}
Let $v=\bbE[(Z-X)^2]=\bbE[X^2]+\bbE[Z^2]$ and $\bbE[X^2]=v_0, \bbE[Z^2]=v_1$, so that $v=v_0+v_1$. Note that if $v=0$, then $X=Z=0$ and it is clear that $E(\rho^*) \le E(\rho)$, and equality holds if only if $\rho=\rho^*$. Hence, from now on we shall assume $v>0$.

Define $f(u,v) := -2u^3-3u^2 -6uv$. Given $v>0$, the function $u \mapsto f(u,v)$ is easily seen to be maximized when $u=\frac{\sqrt{1-4v} -1}{2}$. Plugging in, we obtain
\begin{align} f(u,v) \le g(v) :=-\frac{1}{2} +3v +\frac{1}{2} (1-4v)^{\frac{3}{2}}. \nn
\end{align}
$g(0)=g'(0)=0$ and $g''(0)=6$, hence $g(v) < 3.5v^2 \le 7(v_0^2+v_1^2)$ for all small $v>0$. Note that  $\bbE[(Z-X)^3] = \bbE[Z^3] - \bbE[X^3]$. Now by Proposition \ref{estimate} for the case $n=1$, we have
\begin{align} 7v_0^2+2\bbE[X^3] -  \bbE[|X-X'|^3] \le 0, \q 7v_1^2 -2\bbE[Z^3]  -  \bbE[|Z-Z'|^3] \le 0 \nn
\end{align}
for small $\ep$. This implies $E(\rho^*) < E(\rho)$,  concluding the proof. \qed
\end{pf}

%
%
%
%
%

\section{Further observation}\label{FO}
We can ask the following converse question:
$$\text{If $p$ is small, then does $\rho^*=\frac{1}{2}\delta_0 + \frac{1}{2}\delta_1$ fail to be a minimizer?}$$
For an answer we may consider the following simple competitor
$$\rho_m= \frac{1-m}{2} \delta_0 + m \delta_{\frac{1}{2}}+  \frac{1-m}{2} \delta_1$$
and consider the associated energy
$$E(m):=E(\rho_m)= m(1-m)V(\frac{1}{2}) + \frac{(1-m)^2}{4}V(1).$$
If $E'(0)<0$, then $\rho_0=\rho^*$ is not a $d_\lambda$-local minimizer. Let us see when this should happen. Given $q$, define
$$f(p) := -E'(0)= \frac{V(1)}{2} - V(\frac{1}{2}) = \frac{1}{2p}- \frac{1}{2q} -\frac{1}{p2^p} + \frac{1}{q2^q}.$$
$f(q)=0$, so we ask when $f'(q) >0$. We find
$$f'(q)=\frac{1}{q^22^q}(q\log2 +1 - 2^{q-1}).$$
Let $g(q)=q\log2 +1 - 2^{q-1}$ and note that $g$ is concave. Let $q^* >2$ be the unique positive solution to $g(q)=0$. We summarize as follows.
\begin{prop}\label{converse}
Let $q^*$ be the unique positive solution to the equation $q\log2=2^{q-1}-1$, and let $0<q<q^*$. Then there exists $p_* = p_*(q)> q$ such that  for all $p \in (q, p_*)$ and $1 \le \lambda <\infty$, $\rho^*=\frac{1}{2}\delta_0 + \frac{1}{2}\delta_1$ is not a $d_\lambda$-local minimizer.
\end{prop}
This inspires a few questions, and we leave them for future research.\\
{\bf Q1.} For the pairs $(p,q)$ in Proposition \ref{converse}, what is a global minimizer?\\
{\bf Q2.} If $p>q > q^*$, is  $\rho^*$ a global minimizer? \\
{\bf Q3.} For $q \in [2, q^*)$, is $p_*=p^*$ or not? ($p^*$ is addressed in Remark \ref{remark1}.)\\

 \section*{Acknowledgments}
\noin K. Kang's work is supported by NRF-2017R1A2B4006484 and NRF-
2015R1A5A1009350. H.K. Kim's work is supported by NRF-2018R1D1A1B07049357.  T. Lim gratefully acknowledges support from ShanghaiTech University, and in addition, T. Lim is grateful for the support of the University of Toronto and its Fields Institute for the Mathematical
Sciences, where parts of this work were performed. G. Seo's work is supported by NRF-2017R1A2B4006484. 

\bibliography{joint_biblio}{}

\begin{thebibliography}{4}

\bibitem{ags08} \textsc{L. Ambrosio, N. Gigli,  G. Savar\'{e}},
   \textit{Gradient flows in metric spaces and in the space of probability measures},
Second edition. Lectures in Mathematics ETH Zurich. Birkhauser
Verlag, Basel (2008)



\bibitem{bclr13} \textsc{D. Balagu\'{e}, J. A. Carrillo, T. Laurent , G. Raoul},
   \textit{Dimensionality of Local Minimizers of the Interaction Energy},
Arch. Ration. Mech. Anal. 209(3), 1055-1088 (2013)





\bibitem{bclr2013} \textsc{D. Balagu\'{e}, J. A. Carrillo, T. Laurent , G. Raoul},
   \textit{Nonlocal interactions by repulsive-attractive potentials: radial ins/stability},
Phys. D 260 , 5-25 (2013)



\bibitem{bcp97} \textsc{D. Benedetto, E. Caglioti, M. Pulvirenti  }, \textit{ A kinetic equation for granular media. RAIRO Modél}, Math. Anal. Numér. 31(5), 615-641 (1997)


\bibitem{bcl09} \textsc{A. Bertozzi, J. A. Carrillo, T. Laurent}, \textit{Blow-up in multidimensional aggregation equations with mildly singular interaction kernels}, Nonlinearity 22(3) 683–710 (2009)




\bibitem{bl07} \textsc{A. Bertozzi,  T. Laurent}, \textit{Finite-time blow-up of solutions of an aggregation equation in $\R^n$}, Comm. Math. Phys. 274(3), 717–735 (2007)


\bibitem{bll12} \textsc{A. Bertozzi,  T. Laurent, F. Léger},
   \textit{Aggregation and spreading via the Newtonian potential: the dynamics of patch solutions},
Math. Models Methods Appl. Sci. 22(suppl. 1) 1140005  (2012)

\bibitem{blr11} \textsc{A. Bertozzi, T. Laurent, J. Rosado},
   \textit{$L^p$ theory for the multidimensional aggregation equation},
Comm. Pure Appl. Math. 64(1), 45–83 (2011)













\bibitem{cdfls11} \textsc{J. A. Carrillo,  M. DiFrancesco,  A. Figalli,  T. Laurent,  D. Slep\v{c}ev},
   \textit{Global-in-time weak measure solutions and finite-time aggregation for nonlocal interaction equations},
 Duke Math. J. 156(2), 229-271 (2011)







\bibitem{cfp17} \textsc{J. A. Carrillo
, A. Figalli , F.S. Patacchini},
   \textit{Geometry of minimizers for the interaction energy with mildly repulsive potentials},
   Ann. I. H. Poincar$\acute{e}$ - AN 34  1299-1308 (2017)


\bibitem{ch17}
\textsc{J. A. Carrillo, Y. Huang},
 \textit{Explicit equilibrium solutions for the aggregation equation with power-law potentials}, Kinet. Relat. Models 10(1), 171–192 (2017)


\bibitem{chm14}
\textsc{J. A. Carrillo, Y. Huang, S. Martin},
 \textit{Nonlinear stability of flock solutions in second-order swarming models.}, Nonlinear Anal. Real World Appl. 17, 332–343 (2014)












\bibitem{cmv03} \textsc{J. A. Carrillo,  R. J. McCann, C. Villani},
   \textit{Kinetic equilibration rates for granular media and related equations: entropy dissipation and mass transportation estimates},
Rev. Mat. Iberoamericana 19, 971–1018,  (2003)

\bibitem{cmv06} \textsc{J. A. Carrillo,  R. J. McCann, C. Villani},
   \textit{Contractions in the
2-Wasserstein length space and thermalization of granular media},
Arch. Ration. Mech. Anal. 179(2), 217-263 (2006)



\bibitem{fr10}
\textsc{K. Fellner, G. Raoul},
 \textit{Stable stationary states of non-local interaction equations},
  Math. Models Methods Appl. Sci., 20(12)  2267--2291 (2010)


\bibitem{fh13} \textsc{R. C. Fetecau, Y. Huang},
   \textit{Equilibria of biological aggregations with nonlocal repulsive-attractive interactions},
Phys. D 260, 49–64 (2013)


\bibitem{fhk11} \textsc{R.C. Fetecau, Y. Huang, T. Kolokolnikov},
   \textit{Swarm dynamics and equilibria for a nonlocal aggregation model},
Nonlinearity 24(10), 2681–2716 (2011)


\bibitem{hkp13} \textsc{T. Kolokolnikov, Y. Huang,  M. Pavlovski},
   \textit{Singular patterns for an aggregation model with a confining potential}, Phys. D 260, 65–76 (2013)


   

\bibitem{bksu11} \textsc{T. Kolokolnikov, H. Sun,  D. Uminsky, A. Bertozzi}, \textit{ A theory of complex patterns arising from 2D particle interactions}, Phys. Rev. E, Rapid Commun. 84  015203(R) (2011)



%





\bibitem{me99} \textsc{A. Mogilner, L. Edelstein-Keshet }, \textit{A non-local model for a swarm}, J. Math. Biol. 38(6), 534-570 (1999)



\bibitem{s17} \textsc{F. Santambrogio},
   \textit{{Euclidean, metric, and Wasserstein} gradient flows: an overview},
Bull. Math. Sci. 7(1), 87-154 (2017)



\bibitem{bt04} \textsc{C. Topaz, A. Bertozzi},
   \textit{Swarming patterns in a two-dimensional kinematic model for biological groups},
SIAM J. Appl. Math. 65,  152–174 (2004)


\bibitem{blt06} \textsc{C. Topaz, A. Bertozzi, M. Lewis},
   \textit{A nonlocal continuum model for biological aggregation},
Bull. Math. Biol. 68(7), 1601-1623 (2006)




\bibitem{t00} \textsc{G. Toscani}, \textit{One-dimensional kinetic models of granular flows } M2AN Math. Model. Numer. Anal. 34(6), 1277-1291 (2000)




\bibitem{v03} \textsc{C. Villani},
   \textit{Topics in optimal transportation},
Graduate Studies in Mathematics, 58. American Mathematical Society,
Providence, RI (2003)


\bibitem{bkuv12} \textsc{J. von Brecht, D. Uminsky, T. Kolokolnikov, A. Bertozzi},
   \textit{Predicting pattern formation in particle interactions},
Math. Models Methods Appl. Sci. 22(suppl. 1) 1140002 (2012)
























%
%
%
%
%
%
%
%
%
%
%
%
%
%
%
%
%
%
%
%
%
%
%
%
%
%
%
%
%
%
%
%
%
%
%
%
%
%



%


%
%
%
%
%
%
%
%
%




%
%







%
%
%
%
%
%


















\end{thebibliography}
\bibliographystyle{plain}

\end{document}